\documentclass[10pt]{amsart}
\usepackage{amsmath, amsthm, amssymb}
\usepackage{mathbbol, mathrsfs}
\usepackage{txfonts}

\newtheorem{thm}{Theorem}[section]
\newcommand{\bt}{\begin{thm}}
\newcommand{\et}{\end{thm}}

\newtheorem{cor}[thm]{Corollary}
\newcommand{\bc}{\begin{cor}}
\newcommand{\ec}{\end{cor}}

\newtheorem{lem}[thm]{Lemma}
\newcommand{\bl}{\begin{lem}}
\newcommand{\el}{\end{lem}}

\newtheorem{prop}[thm]{Proposition}
\newcommand{\bp}{\begin{prop}}
\newcommand{\ep}{\end{prop}}

\newtheorem{defn}[thm]{Definition}
\newcommand{\bd}{\begin{defn}}      
\newcommand{\ed}{\end{defn}}

\newtheorem{rmrk}[thm]{Remark}
\newcommand{\br}{\begin{rmrk}}
\newcommand{\er}{\end{rmrk}}

\newtheorem{question}[thm]{Question}
\newcommand{\bq}{\begin{question}}
\newcommand{\eq}{\end{question}}

\newtheorem{example}[thm]{Example}

\newcommand{\thmref}[1]{Theorem~\ref{#1}}
\newcommand{\secref}[1]{Section~\ref{#1}}
\newcommand{\lemref}[1]{Lemma~\ref{#1}}
\newcommand{\defref}[1]{Definition~\ref{#1}}
\newcommand{\corref}[1]{Corollary~\ref{#1}}
\newcommand{\propref}[1]{Proposition~\ref{#1}}

\newcommand{\exref}[1]{Example~\ref{#1}}

\newcommand{\C}{\mathbb{C}}
\newcommand{\N}{\mathbb{N}}

\newcommand{\R}{\mathbb{R}}
\newcommand{\Z}{\mathbb{Z}}
\newcommand{\Qu}{\mathbb{H}}
\newcommand{\K}{\mathbb{K}}
\newcommand{\Li}{\mathbb{L}}
\newcommand{\Oc}{\mathbb{O}}
\newcommand{\Impart}{\operatorname{Im}}
\newcommand{\length}{\operatorname{length}}
\newcommand{\hm}{{\mathcal H}}
\newcommand{\lm}{{\mathcal L}}

\newcommand{\Area}{\operatorname{Area}}

\newcommand{\lip}{\operatorname{Lip}}
\newcommand{\mass}[2][]{{\mathbf M_{#1}}(#2)}

\newcommand{\form}{{\mathcal D}}        %duality with lipschitz functions%
\newcommand{\curr}{{\mathbf M}}         %metric current%
\newcommand{\intcurr}{{\mathbf I}}      %metric integral current%

\newcommand{\area}{\operatorname{Area}}

\newcommand{\fillarea}{{\operatorname{Fillarea}}}

\newcommand{\FA}{{\operatorname{FA}}}

\newcommand{\rstr}{\:\mbox{\rule{0.1ex}{1.2ex}\rule{1.1ex}{0.1ex}}\:}
\newcommand{\bdry}{\partial}

\newcommand{\spt}{\operatorname{spt}}
\newcommand{\ohne}{\backslash}

\begin{document}

\title{Nilpotent groups without exactly polynomial Dehn function}

\author{Stefan Wenger}

\address
{Department of Mathematics\\
University of Illinois at Chicago\\
851 S. Morgan Street\\
Chicago, IL 60607--7045}
\email{wenger@math.uic.edu}

\date{April 16, 2010}

%\keywords{Isoperimetric inequalities, asymptotic rank, Euclidean rank, non-positive curvature, Hadamard spaces, cone type inequalities}

\thanks{Partially supported by NSF grant DMS 0956374}

%49Q15 (Geometric measure and integration theory, integral and normal currents)
%53C23

%\subjclass[2000]{49Q15}

\begin{abstract}
 We prove super-quadratic lower bounds for the growth of the filling area function of a certain class of Carnot groups. This class contains groups for which it is known that their Dehn function grows no faster than  $n^2\log n$. We therefore obtain the existence of (finitely generated) nilpotent groups whose Dehn functions do not have exactly polynomial growth and we thus answer a well-known question about the possible growth rate of Dehn functions of nilpotent groups.
\end{abstract}

\maketitle

\section{Introduction}\label{section:introduction}

Dehn functions have played an important role in geometric group theory. They measure the complexity of the word problem in a given finitely presented group and provide a quasi-isometry invariant of the group. 
A class of groups for which the Dehn function has been well-studied is that of nilpotent groups. 
%Dehn functions of nilpotent groups have been particularly well-studied. 
Many results on upper bounds for the growth of the Dehn function are known, while fewer techniques have been found for obtaining lower bounds. Here are some known results:
if $\Gamma$ is a finitely generated nilpotent group of step $c$ then its Dehn function $\delta_\Gamma(n)$ grows no faster than $n^{c+1}$, in short $\delta_\Gamma(n)\preceq n^{c+1}$, \cite{Gromov-asymptotic, Pittet-homogeneous-nilpotent, Gersten-Holt-Riley}. If $\Gamma$ is a free nilpotent group of step $c$ then $\delta_\Gamma(n)\sim n^{c+1}$, \cite{BMS, Pittet-isop-nilpotent-London}. This includes the particular case of the first Heisenberg group, which was known before, see \cite{Epstein-et-all}. For the definition of $\delta_\Gamma(n)$ and the meaning of $\preceq$ and $\sim$ see \secref{section:prelims}. The higher Heisenberg groups have quadratic Dehn functions, \cite{Gromov-asymptotic, Allcock, Olshanskii-Sapir}. Ol'shanskii-Sapir used the arguments in their paper \cite{Olshanskii-Sapir} to prove that the Dehn function of the central product of $m\geq 2$ copies of a nilpotent group of step $2$ grows no faster than $n^2\log n$. This has also recently been proved in \cite{Young-scaled-relators} with different methods.

The following question about the possible growth rate of the Dehn function of nilpotent groups has been raised by many authors, see e.g.~\cite{Pittet-isop-nilpotent-London, AMS, Riley-filling-notes}, but has remained open so far.

\bq\label{question}
 Does the Dehn function of every finitely generated nilpotent group $\Gamma$ grow exactly polynomially, that is, $$\delta_\Gamma(n)\sim n^\alpha \quad\text{for some (integer) $\alpha$?}$$ 
\eq

The primary purpose of this article is to prove super-quadratic lower bounds for the growth of the Dehn function of certain classes of nilpotent groups and to combine these with the results from \cite{Olshanskii-Sapir, Young-scaled-relators} to give a negative answer to the above question. We establish:

\bt\label{thm:fin-gen-nilpotent}
 There exist finitely generated nilpotent groups $\Gamma$ of step $2$ whose Dehn function $\delta_\Gamma(n)$ satisfies
 \begin{equation}\label{eq:intermediate-growth}
  n^2 \varrho(n) \preceq \delta_\Gamma(n) \preceq n^2\log n
 \end{equation}
 for some function $\varrho$ with $\varrho(n)\to\infty$ as $n\to\infty$.
\et

We can in fact produce a whole family of groups satisfying \eqref{eq:intermediate-growth}. They all arise as lattices of (central powers) of Carnot groups of step $2$, and the upper bound comes as a consequence of \cite{Olshanskii-Sapir, Young-scaled-relators}. Recall that a simply connected nilpotent Lie group is called Carnot group if its Lie algebra admits a grading, the first layer of which generates the whole Lie algebra, see \secref{section:Carnot-groups} for definitions. In \thmref{thm:inaccessible} we will show that under suitable conditions on the Lie algebra of a Carnot group $G$ of step $2$ the $m$-th central power of $G$ has filling area function which grows strictly faster than quadratically.  Recall that the filling area function $\FA(r)$ of a simply connected Riemannian manifold $X$ is the smallest function such that every closed curve in $X$ of length at most $r$ bounds a singular Lipschitz chain of area at most $\FA(r)$. As is well-known, $\FA(r)$ bounds from below the Dehn function of any finitely generated group acting properly discontinuously and cocompactly by isometries on $X$.

In preparation of the proof of the super-quadratic lower bounds and of \thmref{thm:fin-gen-nilpotent}, we will prove, in \secref{section:quad-isop}, the following theorem which should be of independent interest.

%Before describing the conditions which imply super-quadratic lower bounds we wish to elaborate on the general method of proof we develop. The following is a special case of a more general theorem which we prove in 

\bt\label{thm:special-quadratic}
Let $X$ be a simply connected Riemannian manifold with quadratic filling area function, $\FA(r)\preceq r^2$. Then every asymptotic cone $X_\omega$ of $X$ admits a quadratic isoperimetric inequality for integral $1$-currents in $X_\omega$. In particular, there exists $C>0$ such that every closed Lipschitz curve $c$ in $X_\omega$ bounds an integral $2$-current in $X_\omega$ of mass (``area'') at most $C\length(c)^2$.
\et

The above theorem holds more generally for complete metric length spaces $X$, see \thmref{main-thm-intcurr-isop}. The theory of integral currents in metric spaces was developed by Ambrosio-Kirchheim in \cite{Ambr-Kirch-curr}. We refer to \secref{section:prelims} for definitions. Integral $2$-currents should be thought of as generalized singular Lipschitz chains and the mass is proportional to the Hausdorff area. 
A consequence of \thmref{thm:special-quadratic} is the following:

\bc\label{corollary:quad-Dehn-asymp}
 Let $\Gamma$ be a finitely presented group with quadratic Dehn function. Then every asymptotic cone of $\Gamma$ admits a quadratic isoperimetric inequality for integral $1$-currents.
\ec

It is known that finitely presented groups with quadratic Dehn function have simply connected asymptotic cones, \cite{Papasoglu-quadratic}. Many groups, however, have simply connected asymptotic cones but do not have quadratic Dehn function. \corref{corollary:quad-Dehn-asymp} does not yield simple connectedness but instead provides strong metric information. This can be used to prove that certain groups cannot admit quadratic Dehn function, see below.

We now briefly describe how \thmref{thm:special-quadratic} will be used to obtain super-quadratic growth estimates for filling area functions. Let $H$ be a Carnot group and let $d_0$ be the distance induced by a left-invariant Riemannian metric on $H$. Denote by $d_c$ the Carnot-Carath\'eodory distance on $H$ associated with $d_0$, see \secref{section:Carnot-groups}. By \cite{Pansu-boules}, the metric space $(H, d_c)$ is the unique asymptotic cone of $(H, d_0)$. In order to prove that the filling area function of $(H, d_0)$ grows strictly faster than quadratically for a given $(H, d_0)$ it thus suffices, by \thmref{thm:special-quadratic}, to prove that the metric space $(H, d_c)$ cannot admit a quadratic isoperimetric inequality for integral $1$-currents. In Sections \ref{section:inaccessible} and \ref{section-thm-subalgebra} we will show that under suitable conditions on the Lie algebra of $H$, the metric space $(H, d_c)$ has non-trivial first homology group for integral currents, and in particular $(H, d_c)$ does not admit a quadratic isoperimetric inequality for integral $1$-currents. The conditions exhibited in \secref{section:inaccessible} will be used to establish the lower bounds needed to prove \thmref{thm:fin-gen-nilpotent}; they can roughly be described as follows. Let $G$ be a Carnot group of step $2$ with grading $\mathfrak{g}= V_1\oplus V_2$ of its Lie algebra. Given a proper subspace $U\subset V_2$, one obtains a new Lie algebra $\mathfrak{g}_U = V_1\oplus (V_2/U)$ and hence a Carnot group $G_U$ with Lie algebra $\mathfrak{g}_U$. Note that every  simply connected nilpotent Lie group of step $2$ arises this way with $\mathfrak{g}=V_1\oplus V_2$ a free nilpotent Lie algebra of step $2$ and $U\subset V_2$ a suitable subspace.
In the proof of \thmref{thm:inaccessible} we will show that under a suitable condition on $U$, which we term $m$-inaccessibility, the $m$-th central power $H:= G_U\times_Z\dots\times_Z G_U$, endowed with a Carnot-Carath\'eodory distance, has non-trivial first Lipschitz homology. In view of the above, this is enough to conclude that the filling area function of $H$ grows super-quadratically when $H$ is endowed with a left-invariant Riemannian metric.
The $m$-inaccessibility condition is given in \defref{def:inaccessible}; it is not difficult to provide examples of $m$-inaccessible subspaces $U$, see \exref{example:inaccessible} and the remark following it.
For the proofs of Theorems \ref{thm:inaccessible} and \ref{thm:fin-gen-nilpotent}, we only need the special case of \thmref{thm:special-quadratic} when $X$ is a Carnot group endowed with a left-invariant Riemannian metric. In this case, the proof can be simplified and we can in fact prove an analog for cycles of any dimension, see \propref{prop-Carnot-isop}. We believe that, nevertheless, the generality of \thmref{thm:special-quadratic} should be of interest, especially since not much is known about groups with quadratic Dehn function.

The paper is structured as follows: \secref{section:prelims} contains definitions of and facts about Dehn functions, filling area functions, asymptotic cones, and Carnot groups. We will also recall all definitions from Ambrosio-Kirchheim's theory of integral currens \cite{Ambr-Kirch-curr} which we will need later. In  \secref{section:quad-isop} we will use the theory of currents to prove \thmref{thm:special-quadratic} and the more general \thmref{main-thm-intcurr-isop}. In \secref{section:inaccessible}, \thmref{main-thm-intcurr-isop} will be used to prove the super-quadratic lower bounds in \thmref{thm:inaccessible} and to establish \thmref{thm:fin-gen-nilpotent}.
Finally, in \secref{section-thm-subalgebra}, we use \thmref{main-thm-intcurr-isop} to exhibit super-quadratic lower bounds for another class of Carnot groups, not necessarily of step $2$, see \thmref{thm-subalgebra}. The results obtained in \secref{section-thm-subalgebra} are not needed for the proof of \thmref{thm:fin-gen-nilpotent}.

\bigskip

{\bf Acknowledgment:} I am indebted to Robert Young for inspiring discussions and for bringing Question~\ref{question} to my attention. Parts of the research underlying this paper was carried out during a research visit to the Universit\'e Catholique de Louvain, Belgium. I wish to thank the Mathematics Department and Thierry de Pauw for the hospitality.

\section{Preliminaries}\label{section:prelims}

The purpose of this section is to collect some definitions and facts which will be used in the sequel.

\subsection{Dehn and filling area functions}\label{section:prelim-dehn}

Let $\Gamma = \,<S\,|\,R>$ be a finitely presented group. The Dehn function $\delta_\Gamma(n)$ of $\Gamma$ with respect to the given presentation is defined by
\begin{equation*}
 \delta_\Gamma(n):= \max_{w =_\Gamma 1, |w|\leq n} \min\left\{k: w = \prod_{i=1}^k g_i r_i^{\pm 1} g_i^{-1}, r_i\in R, \text{$g_i$ word in $S$}\right\}.
\end{equation*}
Here, $w$ is a word in the alphabet $S$ and $|w|$ denotes the word length of $w$. The equality $w = \prod_{i=1}^k g_i r_i^{\pm 1} g_i^{-1}$ is in the free group generated by $S$.
For functions $f, g:\N\to\N$ or $f, g:[0,\infty]\to (0,\infty)$ one writes $f\preceq g$ if there exists $C$ such that 
\begin{equation*}
 f(r) \leq Cg(Cr+C) + Cr + C
\end{equation*} 
for all $r>0$. One furthermore writes $f\sim g$ if $f\preceq g$ and $g\preceq f$.
If $\Gamma$ acts properly discontinuously and cocompactly by isometries on a simply connected Riemannian manifold $X$ then 
\begin{equation*}
 \delta_\Gamma(n) \sim \FA_0(n),
\end{equation*}
where the Dehn function $\FA_0(r)$ on $X$ is defined by
\begin{equation*}
 \FA_0(r):= \sup\{\fillarea_0(c): \text{$c$ closed Lipschitz curve in $X$ of $\length(c)\leq r$}\},
\end{equation*}
and $\fillarea_0(c)$ of $c:S^1\to X$ is given by
\begin{equation*}
 \fillarea_0(c):= \inf\left\{\area(\varphi): \text{$\varphi: D^2\to X$ Lipschitz, $\varphi|_{S^1} = c$}\right\}.
\end{equation*}
Here, $\Area(\varphi)$ is the integral over the disc $D^2$ of the jacobian of $\varphi$.
The filling area function $\FA(r)$ of $X$ is defined by
\begin{equation*}
 \FA(r):= \sup\{\fillarea(c): \text{$c$ closed Lipschitz curve in $X$ of $\length(c)\leq r$}\},
\end{equation*}
where
\begin{equation*}
 \fillarea(c):= \inf\left\{\Area(z): \text{$z$ is a singular Lipschitz $2$-chain in $X$ with boundary $c$}\right\}.
\end{equation*}
Recall that a singular Lipschitz $2$-chain is a formal finite sum $z = \sum m_i\varphi_i$ with $m_i\in\Z$ and $\varphi_i: \Delta^2\to X$ Lipschitz and its area is
\begin{equation*}
 \Area(z):= \sum |m_i| \Area(\varphi_i),
\end{equation*}
see \cite{Gromov-filling}.
It is clear that $$\FA(r)\leq \FA_0(r)$$ for all $r\geq 0$.

\subsection{Ultralimits and asymptotic cones of metric spaces}

Recall that a non-principal ultrafilter on $\N$ is a finitely additive probability measure $\omega$ on $\N$ (together with the $\sigma$-algebra of all subsets) such
 that $\omega$ takes values in $\{0,1\}$ only and
 $\omega(A)=0$ whenever $A\subset \N$ is finite.
The existence of non-principal ultrafilters on $\N$ follows from Zorn's lemma. It is not difficult to prove the following fact.
If $(Y,\tau)$ is a compact topological Hausdorff space then for every sequence $(y_n)_{n\in\N}\subset Y$ there exists a unique point $y\in Y$ such that
\begin{equation*}
 \omega(\{n\in\N: y_n\in U\})=1
\end{equation*}
for every $U\in\tau$ containing $y$. We will denote this point by $\lim\nolimits_\omega y_n$.\\
Let now $(X_n,d_n, p_n)$ be pointed metric spaces, $n\geq 1$, and fix a non-principal ultrafilter $\omega$ on $\N$. A sequence $(x_n)$, where $x_n\in X_n$ for each $n$, will be called bounded if $$\sup_n d_n(p_n, x_n)<\infty.$$ Define an equivalence relation on bounded sequences by
\begin{equation*}
 (x_n)\sim (x'_n)\quad\text{if and only if}\quad \lim\nolimits_\omega d_n(x_n,x'_n)=0.
\end{equation*}
Let $X_\omega$ be the set of equivalence classes of bounded sequences, $d_\omega$ the metric on $X_\omega$ given by
\begin{equation*}
 d_\omega([(x_n)],[(x'_n)]):= \lim\nolimits_\omega d_n(x_n,x'_n),
\end{equation*}
and $p_\omega:= [(p_n)]$. Then the pointed metric space $(X_\omega, d_\omega, p_\omega)$  is called the ultralimit of the sequence $(X_n,d_n, p_n)$ with respect to $\omega$. If $(X, d)$ is a metric space, $(p_n)\subset X$ a sequence of basepoints, and $r_n>0$ satisfies $r_n\to\infty$, then the ultralimit of the sequence $(X, r_n^{-1}d, p_n)$ with respect to $\omega$ is called the asymptotic cone of $(X, (r_n), (p_n))$ with respect to $\omega$.

It is not difficult to show that ultralimits are always complete. Furthermore, if every $(X_n, d_n)$ is a length space, then the ultralimits are all length spaces as well.

\subsection{Carnot groups and their asymptotic cones}\label{section:Carnot-groups}

Recall that a connected and simply connected nilpotent Lie group $G$ of step $k$ is called Carnot group if its Lie algebra $\mathfrak{g}$ admits a grading $$\mathfrak{g} = V_1\oplus\dots\oplus V_k$$ such that $[V_1, V_i]= V_{i+1}$ for all $i=1, \dots, k-1$ and $[V_1, V_k]=0$, where $[V_1, V_i]$ is the smallest subspace spanned by the elements $[v, v']$ with $v\in V_1$ and $v'\in V_i$. In other words, a Carnot group is a homogeneous nilpotent Lie group such that the first layer $V_1$ of its Lie algebra generates the entire Lie algebra.
Clearly, every connected and simply connected nilpotent Lie group of step $2$ is a Carnot group. 
Let now $G$ be a Carnot group of step $k$ and let $$\mathfrak{g} = V_1 \oplus \dots \oplus V_k$$ be a grading of its Lie algebra $\mathfrak{g}$. Note that the exponential map $\exp: \mathfrak{g}\to G$ is a diffeomorphism. $G$ comes with a family of dilation homomorphisms $\delta_r: G\to G$, $r\geq 0$, which, on the level of Lie algebras, take the form
\begin{equation*}
 \delta_r(v) = \sum_{i=1}^k r^iv_i
\end{equation*}
for $v = v_1+\dots+ v_k$ with $v_i\in V_i$.

Let $g_0$ be a left-invariant Riemannian metric on $G$ and denote by $d_0$ the distance on $G$ coming from $g_0$. A new distance, called Carnot-Carath\'eodory distance, can be associated with $g_0$ as follows. Define the horizontal subbundle $TH$ of $TG$ by left-translating $V_1$. A curve $c:[0,1]\to G$, absolutely continuous (ac for short) with respect to $d_0$, is called horizontal if $\dot{c}(t)\in T_{c(t)}H$ for almost every $t\in[0,1]$. 
The Carnot-Carath\'eodory distance $d_c$ on $G$ associated with $g_0$ is then defined by
\begin{equation*}
 d_c(x,y) = \inf \{ length_{g_0} (c): \text{$c$ horizontal ac curve joining $x$ to $y$}\},
\end{equation*}
where $\length_{g_0}(c)$ denotes the length of $c$ with respect to $g_0$. 
It can be shown that $d_c$ defines a metric, ie. that $d_c(x,y)$ is always finite. Important properties of the Carnot-Carath\'eodory distance are that it is left invariant and 1-homogeneous with respect to the dilations, i.e., $d_c(\delta_r(x),\delta_r(y)) = r\,d_c(x,y)$ for all $x, y \in G$ and all $r\geq 0$.
We also obviously have the relationship
\begin{equation*}
 d_0 \leq d_c,
\end{equation*}
and it is well-known that $d_0$ and $d_c$ are not bi-Lipschitz equivalent unless $G$ is Euclidean. Note however, that the topologies induced by $d_0$ and $d_c$ are the same.

The following theorem, which is a special case of a more general result due to Pansu \cite{Pansu-boules}, gives a link between left-invariant Riemannian metrics and Carnot-Carath\'eodory metrics on  Carnot groups.

\bt\label{pansu-thm}
 Let $G$ be a Carnot group and let $d_0$ be the distance associated with a left-invariant Riemannian metric $g_0$ on $G$. Then the pointed spaces $(G, \frac{1}{r}d_0, e)$ converge in the pointed Gromov-Hausdorff sense to $(G, d_c, e)$ as $r\to\infty$, where $d_c$ is the Carnot-Carath\'eodory distance on $G$ associated with $g_0$. Here, $e$ denotes the identity element of $G$.
\et

In particular, it follows that $(G, d_0)$ has a unique asymptotic cone, which moreover is isometric to $(G, d_c)$. This will be used in the proofs of  \thmref{thm:inaccessible} and  \thmref{thm-subalgebra}.

\subsection{Integral currents in metric spaces}\label{section:currents}
The theory of integral currents in metric spaces was developed by Ambrosio and Kirchheim in \cite{Ambr-Kirch-curr} and provides a suitable 
notion of surfaces and area/volume in the setting of metric spaces. In the following we recall the definitions that are needed for our purposes.

Let $(Y,d)$ be a complete metric space and $m\geq 0$ and let $\form^m(Y)$ be the set of $(m+1)$-tuples $(f,\pi_1,\dots,\pi_m)$ 
of Lipschitz functions on $Y$ with $f$ bounded. The Lipschitz constant of a Lipschitz function $f$ on $Y$ will
be denoted by $\lip(f)$.
\bd
An $m$-dimensional metric current  $T$ on $Y$ is a multi-linear functional on $\form^m(Y)$ satisfying the following
properties:
\begin{enumerate}
 \item If $\pi^j_i$ converges point-wise to $\pi_i$ as $j\to\infty$ and if $\sup_{i,j}\lip(\pi^j_i)<\infty$ then
       \begin{equation*}
         T(f,\pi^j_1,\dots,\pi^j_m) \longrightarrow T(f,\pi_1,\dots,\pi_m).
       \end{equation*}
 \item If $\{y\in Y:f(y)\not=0\}$ is contained in the union $\bigcup_{i=1}^mB_i$ of Borel sets $B_i$ and if $\pi_i$ is constant 
       on $B_i$ then
       \begin{equation*}
         T(f,\pi_1,\dots,\pi_m)=0.
       \end{equation*}
 \item There exists a finite Borel measure $\mu$ on $Y$ such that
       \begin{equation}\label{equation:mass-def}
        |T(f,\pi_1,\dots,\pi_m)|\leq \prod_{i=1}^m\lip(\pi_i)\int_Y|f|d\mu
       \end{equation}
       for all $(f,\pi_1,\dots,\pi_m)\in\form^m(Y)$.
\end{enumerate}
\ed
The space of $m$-dimensional metric currents on $Y$ is denoted by $\curr_m(Y)$ and the minimal Borel measure $\mu$
satisfying \eqref{equation:mass-def} is called mass of $T$ and written as $\|T\|$. We also call mass of $T$ the number $\|T\|(Y)$ 
which we denote by $\mass{T}$.
The support of $T$ is, by definition, the closed set $\spt T$ of points $y\in Y$ such that $\|T\|(B(y,r))>0$ for all $r>0$. 

Every function $\theta\in L^1(K,\R)$ with $K\subset\R^m$ Borel measurable induces an element of $\curr_m(\R^m)$ by
\begin{equation*}
 \Lbrack\theta\Rbrack(f,\pi_1,\dots,\pi_m):=\int_K\theta f\det\left(\frac{\partial\pi_i}{\partial x_j}\right)\,d\lm^m
\end{equation*}
for all $(f,\pi_1,\dots,\pi_m)\in\form^m(\R^m)$.

The restriction of $T\in\curr_m(Y)$ to a Borel set $A\subset Y$ is given by 
\begin{equation*}
  (T\rstr A)(f,\pi_1,\dots,\pi_m):= T(f\chi_A,\pi_1,\dots,\pi_m).
\end{equation*}
This expression is well-defined since $T$ can be extended to a functional on tuples for which the first argument lies in 
$L^\infty(Y,\|T\|)$.

If $m\geq 1$ and $T\in\curr_m(Y)$ then the boundary of $T$ is the functional
\begin{equation*}
 \bdry T(f,\pi_1,\dots,\pi_{m-1}):= T(1,f,\pi_1,\dots,\pi_{m-1}).
\end{equation*}
It is clear that $\bdry T$ satisfies conditions (i) and (ii) in the above definition. If $\bdry T$ also satisfies (iii) then $T$ is called a normal current.
By convention, elements of $\curr_0(Y)$ are also called normal currents.

The push-forward of $T\in\curr_m(Y)$ 
under a Lipschitz map $\varphi$ from $Y$ to another complete metric space $Z$ is given by
\begin{equation*}
 \varphi_\# T(g,\tau_1,\dots,\tau_m):= T(g\circ\varphi, \tau_1\circ\varphi,\dots,\tau_m\circ\varphi)
\end{equation*}
for $(g,\tau_1,\dots,\tau_m)\in\form^m(Z)$. This defines a $m$-dimensional current on $Z$.
It follows directly from the definitions that $\bdry(\varphi_{\#}T) = \varphi_{\#}(\bdry T)$.

\sloppy
We will mainly be concerned with integral currents. Let $\hm^m$ denote Hausdorff $m$-measure on $Y$ and recall that an $\hm^m$-measurable set $A\subset Y$
is said to be countably $\hm^m$-rectifiable if there exist countably many Lipschitz maps $\varphi_i :B_i\longrightarrow Y$ from subsets
$B_i\subset \R^m$ such that
\begin{equation*}
\hm^m\left(A\ohne \bigcup \varphi_i(B_i)\right)=0.
\end{equation*}
\fussy

An element $T\in\curr_0(Y)$ is called integer rectifiable if there exist finitely many points $y_1,\dots,y_n\in Y$ and $\theta_1,\dots,\theta_n\in\Z\ohne\{0\}$ such
that
\begin{equation*}
 T(f)=\sum_{i=1}^n\theta_if(y_i)
\end{equation*}
for all bounded Lipschitz functions $f$.
 A current $T\in\curr_m(Y)$ with $m\geq 1$ is said to be integer rectifiable if the following properties hold:
 \begin{enumerate}
  \item $\|T\|$ is concentrated on a countably $\hm^m$-rectifiable set and vanishes on $\hm^m$-neg\-li\-gible Borel sets.
  \item For any Lipschitz map $\varphi:Y\to\R^m$ and any open set $U\subset Y$ there exists $\theta\in L^1(\R^m,\Z)$ such that 
    $\varphi_\#(T\rstr U)=\Lbrack\theta\Rbrack$.
 \end{enumerate}
Integer rectifiable normal currents are called integral currents. The corresponding space is denoted by $\intcurr_m(Y)$. If $A\subset\R^m$ is a Borel set of finite measure and
finite perimeter then $\Lbrack\chi_A\Rbrack \in\intcurr_m(\R^m)$. Here, $\chi_A$ denotes the characteristic function. If $T\in\intcurr_m(Y)$ and if $\varphi:Y\to Z$ is a Lipschitz 
map into another complete metric space then $\varphi_{\#}T\in\intcurr_m(Z)$.

\medskip

We close this section with a few remarks. A Lipschitz curve $c:[a,b]\to Y$ gives rise to the element $c_{\#}\Lbrack\chi_{[a,b]}\Rbrack\in\intcurr_1(Y)$.
If $\varphi$ is one-to-one then 
$\mass{c_{\#}\Lbrack\chi_{[a,b]}\Rbrack}=\length(c)$. It was shown in \cite[Lemma 2.3]{Wenger-Gromov-hyp-isop} that $1$-dimensional integral currents are essentially induced by (countably many) Lipschitz curves. 
A Lipschitz map $\varphi:D^2\to Y$ gives rise to the $2$-dimensional integral current $S:= \varphi_{\#}\Lbrack\chi_{D^2}\Rbrack$. If $\varphi$ is one-to-one then $$\frac{1}{c}\Area(\varphi)\leq \mass{S}\leq c\Area(\varphi)$$ for some universal $c$;  if $Y$ is a Riemannian manifold then  $\mass{S} = \Area(\varphi)$.
A singular Lipschitz chain $c=\sum m_i\varphi_i$ gives rise to the integral current $\sum m_i\varphi_{i\#}\Lbrack\chi_{\Delta}\Rbrack$. An element $T\in\intcurr_2(Y)$ can be thought of as a generalized singular Lipschitz chain whose boundary consists of a union of closed Lipschitz curves of finite total length.

\section{Quadratic isoperimetric inequalities and asymptotic cones}\label{section:quad-isop}

A complete metric space $X$ will be said to admit a quadratic isoperimetric inequality for $\intcurr_1(X)$ if 
there exists $C>0$ such that for every $T\in\intcurr_1(X)$ with $\bdry T=0$ there exists $S\in\intcurr_2(X)$ with
 \begin{equation*}
  \mass{S}\leq C \mass{T}^2.
 \end{equation*}
The number $C$ will be called isoperimetric constant.
The main result of this section is the following theorem which will be needed in the proof of \thmref{thm:fin-gen-nilpotent} and which should be of independent interest.

\bt\label{main-thm-intcurr-isop}
 Let $X$ be a complete metric length space. If $X$ admits a quadratic isoperimetric inequality for $\intcurr_1(X)$ with isoperimetric constant $C$ then every asymptotic cone $X_\omega$ of $X$ admits a quadratic isoperimetric inequality for $\intcurr_1(X_\omega)$ with isoperimetric constant $4C$.
\et

For the proof of \thmref{thm:fin-gen-nilpotent} we will actually only need \thmref{main-thm-intcurr-isop} for $X$ a Carnot group endowed with a left-invariant Riemannian metric.
In this case, the proof can be simplified, see \propref{prop-Carnot-isop}. We believe that \thmref{main-thm-intcurr-isop} should be of independent interest.
As a consequence of \thmref{main-thm-intcurr-isop} we obtain the following result.

\bc\label{cor:coarse-quad-asymp}
 Let $X$ be a metric length space which admits a coarse homological isoperimetric inequality for curves. Then every asymptotic cone $X_\omega$ of $X$ admits a quadratic isoperimetric inequality for $\intcurr_1(X_\omega)$.
 \ec

 The definition of coarse homological quadratic isoperimetric inequality for curves is given in \cite{Wenger-Gromov-hyp-isop}. It is a homological analog of the notion introduced in \cite[III.H.2]{Bridson-Haefliger}.
 
 \begin{proof}[Proof of \corref{cor:coarse-quad-asymp}]
By Proposition 3.2  and Lemma 2.3 of \cite{Wenger-Gromov-hyp-isop}, $X$ has a thickening $X_\delta$ which is a complete metric length space admitting a quadratic isoperimetric inequality for  $\intcurr_1(X_\delta)$. By definition, a thickening of $X$ is a metric space which contains $X$ isometrically and which is at finite Hausdorff distance from $X$. Now, since asymptotic cones of $X_\delta$ and $X$ are the same, \thmref{main-thm-intcurr-isop} shows that every asymptotic cone of $X$ has a quadratic isoperimetric inequality for integral $1$-currents.
\end{proof}

 A special case of the above corollary is the following result.

\bc
 Let $\Gamma$ be a finitely presented group with quadratic Dehn function. Then every asymptotic cone of $\Gamma$ admits a quadratic isoperimetric inequality for integral $1$-currents.
\ec

It is known that every asymptotic cone of a finitely presented group with quadratic Dehn function is simply connected, \cite{Papasoglu-quadratic}. However, many groups have simply connected asymptotic cones even though their Dehn function is not quadratic; for example, the asymptotic cone of any Carnot group is simply connected. \thmref{main-thm-intcurr-isop} does not yield simple connectedness but instead gives strong metric information about the asymptotic cones.

We turn to the proof of \thmref{main-thm-intcurr-isop} for which we need two simple lemmas.

\bl\label{lemma-reduction}
 Let $Y$ be a complete metric length space and $D>0$ and suppose $Y$ satisfies the following condition: for every closed Lipschitz curve $c: [0,1]\to Y$ and every finite partition $0=t_1<t_2<\dots < t_m=1$ of $[0,1]$ there exist a Lipschitz curve $c':[0,1]\to Y$ and $S\in\intcurr_2(Y)$ satisfying $c'(t_i) = c(t_i)$ and
 \begin{equation*}
   \length\left(c'|_{[t_i, t_{i+1}]}\right) \leq  \length\left(c|_{[t_i, t_{i+1}]}\right)
 \end{equation*}
 for all $i$, as well as $\bdry S = c'_\#\Lbrack \chi_{[0,1]}\Rbrack$ and
 \begin{equation*}
  \mass{S}\leq D \length(c)^2.
 \end{equation*}
Then $Y$ admits a quadratic isoperimetric inequality for $\intcurr_1(Y)$ with isoperimetric constant $2D$.
\el

\begin{proof}
 Let $c:[0,1]\to Y$ be a closed Lipschitz curve and set $T:= c_\#\Lbrack\chi_{[0,1]}\Rbrack$. We will find, for each $n\geq 0$, closed Lipschitz curves $c_{n, i}:[0, 1]\to Y$, $i=1, \dots, 8^n$, and $S_n\in\intcurr_2(Y)$ such that 
 $$\length(c_{n,i})\leq \frac{1}{4^n}\length(c)\quad\text{ and }\quad \mass{S_n} \leq \frac{1}{2^{n-1}}D\length(c)^2,$$
as well as
 $$T = \bdry (S_0 + \dots + S_n) + \sum_{i=1}^{8^n} c_{n, i\,\#}\Lbrack\chi_{[0,1]}\Rbrack.$$
For $n=0$ we simply set $c_{0,1}:= c$ and $S_0:= 0$. Suppose we have found $S_0, \dots, S_{n-1}$ and $c_{n-1, i}$ for some $n\geq 1$ with the properties listed above. In order to construct $c_{n, j}$ and $S_n$, we proceed as follows. First, fix $i:=1$ and choose $0=t_1<t_2<\dots < t_9=1$ such that 
$$\length(c_{n-1, i}|_{[t_j, t_{j+1}]}) = \frac{1}{8} \length(c_{n-1, i})$$ for $j=1, \dots, 8$. By assumption, there exist a closed Lipschitz curve $c':[0,1]\to Y$ and $S\in\intcurr_2(Y)$ such that $c'(t_j) = c_{n-1, i}(t_j)$ and
 \begin{equation*}
   \length\left(c'|_{[t_j, t_{j+1}]}\right) \leq  \length\left(c_{n-1, i}|_{[t_j, t_{j+1}]}\right)
 \end{equation*}
 for all $j$, as well as $\bdry S = c'_\#\Lbrack\chi_{[0,1]}\Rbrack$ and
 \begin{equation*}
  \mass{S}\leq D \length(c_{n-1, i})^2.
 \end{equation*}
 Set $S_{n, 1}:= S$ and, for $j=1, \dots, 8$, let $c_{n, j}$ be the concatenation of $c_{n-1, i}|_{[t_j, t_{j+1}]}$ with $(c'|_{[t_j, t_{j+1}]})^-$. It is clear that
 \begin{equation}\label{sums-curves}
  c_{n-1, 1\,\#}\Lbrack\chi_{[0,1]}\Rbrack = \sum_{j=1}^8 c_{n, j\,\#}\Lbrack\chi_{[0,1]}\Rbrack + \bdry S_{n,1}.
 \end{equation}
 Now, do the same for $i=2, \dots, 8^{n-1}$ in order to obtain, after relabeling indices, curves $c_{n, j}$, $j=1, \dots, 8^{n}$, and currents $S_{n,j}$, $j=1, \dots, 8^{n-1}$. Clearly, we have $\length(c_{n, j}) \leq 4^{-n} \length(c)$ and, for $$S_n:= S_{n, 1} + \dots + S_{n, 8^{n-1}},$$ we have
 \begin{equation*}
  \mass{S_n} \leq D\sum_{i=1}^{8^{n-1}} \length(c_{n-1, i})^2 \leq 8^{n-1}D\cdot \left(\frac{1}{4^{n-1}} \length(c)\right)^2 = \frac{1}{2^{n-1}} D\length(c)^2.
 \end{equation*}
Finally, the fact that
\begin{equation*}
   T = \bdry (S_0 + \dots + S_n) + \sum_{i=1}^{8^{n}} c_{n, i\,\#}\Lbrack\chi_{[0,1]}\Rbrack
\end{equation*}
follows from \eqref{sums-curves}. This proves the existence of $c_{n,i}$ and $S_n$ with the properties stated above. Now, set $S:= \sum_{n=1}^\infty S_n$ and note that $S$ is an integer rectifiable $2$-current satisfying
\begin{equation}\label{eq:mass-bd-isop}
 \mass{S} \leq \sum_{n=1}^\infty \mass{S_n} \leq 2D\length(c)^2.
\end{equation}
It is not difficult to show that $\bdry S = T$. Indeed, view $Y$ as a subset of $Y':= l^\infty(Y)$. It is well-known, see \cite{Wenger-GAFA}, that $Y'$ has a quadratic isoperimetric inequality for $\intcurr_1(Y')$. Let $D'$ be the isoperimetric constant. For each $n$ and $i=1, \dots, 8^n$, there thus exists $Q_{n, i}\in\intcurr_2(Y')$ with $\bdry Q_{n,i} = c_{n, i\,\#}\Lbrack\chi_{[0,1]}\Rbrack$ and such that $\mass{Q_{n,i}}\leq D'\length(c_{n,i})^2$. Set $Q_n:= Q_{n,1} + \dots + Q_{n, 8^n}$, note that $Q_n\in\intcurr_2(Y')$, \begin{equation*}
 \mass{Q_n} \leq D'\sum_{i=1}^{8^n}\length(c_{n, i})^2 \leq \frac{1}{2^n} D'\length(c)^2,
\end{equation*}
and $T = \bdry (S_1 + \dots + S_n) + \bdry Q_n$. Since $Q_n$ converges to $0$ in mass and $S_1+\dots+S_n$ converges in mass to $S$, their boundaries converge weakly to $0$ and $\bdry S$, respectively. This shows that indeed $\bdry S = T$. Since $S$ satisfies \eqref{eq:mass-bd-isop} and since $c$ was arbitrary, Lemma 2.3 of \cite{Wenger-Gromov-hyp-isop} shows that $Y$ admits a quadratic isoperimetric inequality for $\intcurr_1(Y)$ with isoperimetric constant $2D$.
\end{proof}

\bl\label{lemma-Hausdorff-limit-ultralimit}
 Let $Z$ be a compact metric space, $m\in\N$, and $\omega$ a non-principal ultrafilter on $\N$. Suppose $A_n\subset Z$ are closed subsets and $a^1_n,\dots, a^m_n\in A_n$ for $n\geq 1$. Then there exists a subsequence $(A_{n_j})$ such that $A_{n_j}$ converges in the Hausdorff distance to a closed subset of $$A:= \{\lim\nolimits_\omega a_n: \text{$a_n\in A_n$ for every $n$}\}$$ and
$a^i_{n_j} \to \lim\nolimits_\omega a^i_n$ as $j\to\infty$, for all $i=1, \dots, m$.
\el

\begin{proof}
 Since $Z$ is compact, there exist integers $m\leq m_1<m_2<\dots$ and, for each $n\in\N$, a sequence $(a_n^j)_{j\geq m+1}\subset A_n$ of points such that $\{a_n^j: j=1, \dots, m_i\}$ is $2^{-i}$-dense in $A_n$. For each $i$ set $a^i:= \lim_\omega a_n^i$, and denote by $C$ the closure of $\{a^i: i\in\N\}$; clearly $C\subset A$. Define for each $j\in\N$
\begin{equation*}
 \Omega_j:= \{n: d(a_n^i, a^i)\leq 2^{-j}\text{ for $i=1, \dots, m_j$}\}
\end{equation*}
and note that $\Omega_1\supset \Omega_2\supset\dots$ and that $\omega(\Omega_j)=1$ for every $j$; in particular, $\Omega_j$ is not finite. Choose $n_1<n_2<\dots$ with $n_j\in\Omega_j$ for all $j$. It follows that $A_{n_j}$ converges in the Hausdorff sense to $C$. Furthermore, we have that $a_{n_j}^i$ converges to $a^i$ as $j\to\infty$ for each $i$ and, in particular, for $i=1,\dots, m$. This concludes the proof.
\end{proof}

We now use the lemmas above to prove the main theorem of this section.

\begin{proof}[{Proof of \thmref{main-thm-intcurr-isop}}]
 Let $\omega$ be a non-principal ultrafilter on $\N$, let $p_n\in X$, and $r_n>0$ with $r_n\to\infty$. Let $X_\omega$ be the asymptotic cone associated with the pointed sequence $(X, r^{-1}_nd, p_n)$ and $\omega$. Denote by $X_n$ the space $X$ endowed with the metric $d_n:= r_n^{-1}d$.
 Let $c:[0,1]\to X_\omega$ be a closed Lipschitz curve and $0=t_1<t_2<\dots< t_m=1$ a partition of $[0,1]$. Set $x^i:= c(t_i)$. Then $x^i = [(x_n^i)]$ for some $x_n^i\in X$ with
 \begin{equation*}
  \sup_n d_n(x_n^i, p_n) <\infty,
 \end{equation*}
and we may assume that $x_n^m = x_n^1$ for all $n$.
 Let $c_n:[0, 1]\to X_n$ be a Lipschitz curve such that for all $i$ we have $c_n(t_i) = x_n^i$ and $c_n|_{[t_i, t_{i+1}]}$ is parametrized proportional to arc-length with $$\length\left(c_n|_{[t_i, t_{i+1}]}\right)\leq d_n(x_n^i, x_n^{i+1}) + \frac{1}{mn}.$$
Note that
 \begin{equation*}
  \sup_n \lip(c_n) < \infty
 \end{equation*}
 and that there exists $R>0$ such that $c_n$ has image in the ball $B_{X_n}(p_n, R)$ for all $n$. It follows in particular that $\sup_n \length(c_n)<\infty$. Set $T_n:= c_{n\#}\Lbrack\chi_{[0,1]}\Rbrack$ and note that $T_n\in\intcurr_1(X_n)$ satisfies $\bdry T_n = 0$ and $\mass{T_n}\leq \length(c_n) + 1/n$. Set $D:=2C$ where $C$ is the isoperimetric constant for $X$. By Lemma 3.4 in \cite{Wenger-GAFA}, see also Theorem 10.6 in \cite{Ambr-Kirch-curr},  there exists $S_n\in\intcurr_2(X_n)$ such that $\bdry S_n = T_n$ and $$\mass{S_n}\leq D\mass{T_n}^2$$ and such that the sequence $(Y_n)$ of metric spaces $Y_n$ given by $$Y_n:= (\spt S_n \cup c_n([0,1])\cup \{p_n\}, d_n)$$ is uniformly compact in the sense of Gromov. It thus follows from Gromov's compactness theorem \cite{Gromov-poly-growth} that there exist a compact metric space $(Z, d_Z)$ and isometric embeddings $\varphi_n: Y_n\to Z$ for all $n$. Set $A_n:= \varphi_n(Y_n)$ and $a^i_n:= \varphi_n(c_n(t_i))$ for $i=1, \dots, m$. By \lemref{lemma-Hausdorff-limit-ultralimit} there exists a subsequence $(A_{n_j})$ such that $A_{n_j}$ converges in the Hausdorff sense to a closed subset of $A:= \{\lim\nolimits_\omega a_n: \text{$a_n\in A_n$ for all $n$}\}$ and $a^i_{n_j}$ converges to $a^i:= \lim_\omega a^i_n$ for every $i$. After possibly passing to a further subsequence, we may assume that $\psi_j:= \varphi_{n_j}\circ c_{n_j}$ converges uniformly to a Lipschitz curve $\psi: [0,1]\to Z$ and that $\varphi_{n_j\#}S_{n_j}$ converges weakly to some $S'\in\intcurr_2(Z)$. The second assertion is a consequence of the closure and compactness theorems for integral currents in compact metric spaces, see Theorems 5.2 and 8.5 in \cite{Ambr-Kirch-curr}. Note also that $\bdry S' = \psi_{\#}\Lbrack\chi_{[0,1]}\Rbrack$; furthermore $\psi(t_i) = a^i$ and
\begin{equation*}
 \length(\psi_j) \leq \frac{1}{n_j} + \sum_{i=1}^{m-1} d_Z(a^i_{n_j}, a^{i+1}_{n_j}).
\end{equation*} 
Since $a^i_{n_j}\to a^i$ and $d_Z(a^i, a^{i+1}) = d_\omega(x^i, x^{i+1})$, we conclude
\begin{equation*}
 \limsup_{j\to\infty} \length(\psi_j) \leq \sum_{i=1}^{m-1} d_\omega(x^i, x^{i+1}) \leq\length(c),
\end{equation*}
hence also $\length(\psi)\leq \length(c)$ and
\begin{equation*}
\mass{S'} \leq D \liminf_{n\to\infty}\mass{\varphi_{n_j\#}T_n}^2\leq \liminf_{n\to\infty} \length(\psi_j)^2
  \leq D\length(c)^2.
\end{equation*}
Since $A_{n_j}$ converges in the Hausdorff sense to a closed subset of $A$ we furthermore have $$\psi([0,1])\cup \spt S' \subset A.$$ 
 Let now $\Lambda: A \to X_\omega$ be the isometric embedding defined by $\Lambda(a):= [(y_n)]$, where $(y_n)\subset X$ is a sequence with $y_n\in Y_n$ and $a = \lim_\omega\varphi_n(y_n)$. It is not difficult to show that $\Lambda$ is well-defined and is an isometric embedding. Furthermore, the Lipschitz curve $c':= \Lambda\circ\psi$ satisfies $c'(t_i) = x^i= c(t_i)$ and $$\length(c')\leq\length(c).$$ 
 We set $S:= \Lambda_{\#}S'$ and note that $\bdry S = c'_{\#}\Lbrack\chi_{[0,1]}\Rbrack$ and $\mass{S}\leq D\length(c)^2$. In view of \lemref{lemma-reduction} the proof is complete.
\end{proof}

As mentioned above, in the special case when $X$ is a Carnot group endowed with a left-invariant Riemannian metric, the proof of \thmref{main-thm-intcurr-isop} can be simplified. In fact, we can prove an analog for higher-dimensional cycles.
Let $m\geq 1$. A complete metric space $X$ is said to admit an isoperimetric inequality of Euclidean type for $\intcurr_m(X)$ if 
there exists $C>0$ such that for every $T\in\intcurr_m(X)$ with $\bdry T=0$ there exists $S\in\intcurr_{m+1}(X)$ with
 \begin{equation*}
  \mass{S}\leq C \mass{T}^{\frac{m+1}{m}}.
 \end{equation*}
If the above holds only for $T$ with $\spt T$ compact, then $X$ will be said to admit an isoperimetric inequality of Euclidean type for compactly supported integral $m$-currents.

 Let $G$ be a Carnot group and $d_0$ the metric on $G$ coming from a left-invariant Riemannian metric, and let $d_c$ be the associated Carnot-Carath\'eodory metric. 

\bp\label{prop-Carnot-isop}
Let $m\geq 1$. If $X:= (G, d_0)$ admits an isoperimetric inequality of Euclidean type for $\intcurr_m(X)$ then $Y:= (G, d_c)$ admits an isoperimetric inequality of Euclidean type for compactly supported integral $m$-currents.
\ep

Note that in the case $m=1$, \propref{prop-Carnot-isop} together with \cite[Lemma 2.3]{Wenger-Gromov-hyp-isop} yield a quadratic isoperimetric inequality for $\intcurr_1(Y)$, that is, also for integral $1$-currents whose supports are not compact. 

\begin{proof}
Let $T\in\intcurr_m(Y)$ with $\spt T$ compact and $\bdry T= 0$. Denote by $\varphi: Y\to X$ the identity map and note that $\varphi$ is $1$-Lipschitz. For each $n\geq 1$ define $T_n:= (\varphi\circ \delta_n)_\# T\in\intcurr_m(X)$. Clearly, we have $\bdry T_n=0$ and $\mass{T_n}\leq n^m\mass{T}. $ By Lemma 3.4 in \cite{Wenger-GAFA} (see also Theorem 10.6 in \cite{Ambr-Kirch-curr}), there exists $S_n\in \intcurr_{m+1}(X)$ such that $\bdry S_n = T_n$,
 \begin{equation*}
  \mass{S_n}\leq C \mass{T_n}^{\frac{m+1}{m}}\leq C n^{m+1} \mass{T}^{\frac{m+1}{m}},
 \end{equation*} 
 and
 \begin{equation*}
  \|S_n\| (B(x, r)) \geq C'r^{m+1}
 \end{equation*}
 for all $x\in\spt S_n$ and $0\leq r\leq d(x, \spt T_n)$. Here, $C$ and $C'$ are constants only depending on the isoperimetric constant for $\intcurr_m(X)$. It follows that there exists $L$ such that $\spt S_n\subset B(e, nL)$ for every $n$, where $e$ is the identity in $G$. Define metric spaces by 
 \begin{equation*}
  Y_n:= (\spt S_n, \frac{1}{n}d_0).
 \end{equation*}
 It follows that $(Y_n)$ is uniformly compact and thus, by Gromov's compactness theorem, there exists a compact metric space $Z$ and isometric embeddings $\psi_n: Y_n\to Z$ for every $n$. Define subsets $A_n:= \psi_n(\spt S_n)$ of $Z$ and maps $\varrho_n:= \psi_n\circ\varphi\circ\delta_n : (\spt T, d_c)\to Z$. After possibly passing to a subsequence, we may assume that $(A_n)$ converges to a closed subset $A\subset Z$ in the Hausdorff sense and $\varrho_n$ converges uniformly to a $1$-Lipschitz map $\varrho: (\spt T, d_c)\to Z$. After possibly passing to a further subsequence, we may assume by the compactness and closure theorems (see Theorems 5.2 and 8.5 in \cite{Ambr-Kirch-curr}) that $\psi_{n\#}S_n$ converges weakly to some $\hat{S}\in\intcurr_{m+1}(Z)$. Let $\omega$ be a non-principal ultrafilter on $\N$ and define a map $\eta: A\to Y$ as follows. Given $a\in A$, let $x_n\in \spt S_n$ such that $\psi_n(x_n)\to a$, and define $$\eta(a):= \lim\nolimits_\omega \delta_{\frac{1}{n}}(x_n).$$ It is not difficult to show that $\eta$ is well-defined, an isometric embedding, and satisfies $\eta\circ\varrho = \operatorname{id}_{\spt T}$. Set $S:= \eta_{\#}\hat{S}$ and note that $\bdry S = T$ as well as $$\mass{S}\leq \liminf_{n\to\infty} \frac{1}{n^{m+1}}\mass{S_n} \leq C\mass{T}^{\frac{m+1}{m}}.$$ This concludes the proof.
\end{proof}

\section{A lower bound for the filling area function of Carnot groups of step $2$}\label{section:inaccessible}

In this section we prove new lower bounds for the filling area function of the central product of certain Carnot groups of step $2$. The main result is \thmref{thm:inaccessible}, which we will use to prove \thmref{thm:fin-gen-nilpotent}.

Let $G$ be a Carnot group of step $2$ with grading $\mathfrak{g} = V_1\oplus V_2$ of its Lie algebra $\mathfrak{g}$ and Lie bracket $[\cdot,\cdot]$. If $U$ is a subspace of $V_2 $ then $[\cdot, \cdot]$ naturally induces a Lie bracket $[\cdot, \cdot]_U$ on $\mathfrak{g}_U = V_1\oplus V'_2$, where $V'_2:= V_2/U$ is the quotient space. Let $G_U$ be the Carnot group whose Lie algebra is $\mathfrak{g}_U$. Note that every connected and simply connected nilpotent Lie group of step $2$ is of the form $G_U$ with $\mathfrak{g}=V_1\oplus V_2$ a free nilpotent Lie algebra of step $2$ and $U\subset V_2$ a suitable subspace.
We will use the following terminology.

\bd\label{def:inaccessible}
Let $m\geq 1$. A non-trivial subspace $U$ of $V_2$ is called $m$-inaccessible if there exists a proper subspace $U'$ of $U$, possibly $U'=\{0\}$, such that if $v_1, \dots, v_m, w_1,\dots, w_m\in V_1$ then the vector $ [v_1, w_1] + \dots + [v_m, w_m]$ is contained in $U$ if and only it is contained in $U'$. 
\ed

It is not difficult to give examples of graded nilpotent Lie algebras $\mathfrak{g} = V_1\oplus V_2$ such that $V_2$ has a non-trivial $m$-inaccessible subspace. 

\begin{example}\label{example:inaccessible}{\rm
Let $k, m$ satisfy $k\geq 2(m+1)$ and let $\mathfrak{g}=V_1\oplus V_2$ be the free nilpotent Lie algebra of step $2$ with $\dim V_1=k$. Then there exist a basis $\{e_1, \dots, e_k\}$ of $V_1$ and a basis $\{e_{i,j}: 1\leq i < j\leq k\}$ of $V_2$ such that the Lie bracket on $\mathfrak{g}$ satisfies $[e_i, e_j] = e_{i,j}$ whenever $i<j$. It is straightforward to check that the one-dimensional subspace $$U:= \operatorname{span}\{e_{1,2} + e_{3,4}+\dots + e_{2m+1,2m+2}\}$$ of $V_2$ is $m$-inaccessible. One may in fact take $U'=\{0\}$. Note that here $\mathfrak{g}_U$ has a basis with rational structure constants.
}
\end{example}

More generally, if $\mathfrak{g} = V_1\oplus V_2$ is a stratified nilpotent Lie algebra of step $2$ with 
\begin{equation}\label{dim-cond}
\dim V_2 > \left(2\dim V_1 - 1\right)m
\end{equation}
then $V_2$ possesses an $m$-inaccessible subspace $U$; furthermore, $U$ may be chosen such that $\mathfrak{g}_U$ has a basis with rational structure constants. Indeed, the subset $C_m\subset V_2$ given by 
\begin{equation*}
 C_m:= \{[v_1, w_1] + \dots + [v_m, w_m]: v_i, w_i\in V_1\}
\end{equation*}
is a cone containing $0$ and is the image of a smooth map $\psi: (\R\times S^k\times S^k)^{m}\to V_2$, where $k=\dim V_1 - 1$ . Consequently,  $C_m$ has Hausdorff dimension at most $(2\dim V_1 - 1)m$ and hence $C_m\not=V_2$ as soon as \eqref{dim-cond} holds. Clearly, for $v\in V_2\backslash C_m$, the subspace $U:= \operatorname{span}\{v\}$ is $m$-inaccessible.

We turn to the main result of this section.

\bt\label{thm:inaccessible}
 Let $G$ be a Carnot group of step $2$ with grading $\mathfrak{g} = V_1\oplus V_2$ of its Lie algebra. Suppose $U\subset V_2$ is an $m$-inaccessible subspace where $m\geq 1$. 
Let $G_U$ be the Carnot group of step $2$ whose Lie algebra is $\mathfrak{g}_U= V_1\oplus V'_2$, where $V'_2:= V_2/U$, and let $$H:= G_U\times_Z\dots\times_Z G_U$$ be the central product of $m$ copies of $G_U$. Then, endowed with a left-invariant Riemannian metric, $H$ has filling area function which grows strictly faster than quadratically:
 \begin{equation*}
  \frac{\FA(r)}{r^2}\to\infty\quad\text{as $r\to\infty$.}
 \end{equation*}
\et

For the definition of $\FA(r)$ see \secref{section:prelim-dehn}.
Recall that the central product of $m$ copies of a group $\Gamma$ is the quotient
of the $m$-fold direct product of $\Gamma$ by the normal subgroup $N$ of tuples $(g_1, \dots, g_m)$ with $g_i\in[\Gamma, \Gamma]$ and $g_1\cdots g_m = e$, where $e$ is the identity element of $\Gamma$. In order to prove \thmref{thm:inaccessible} we will actually show that there exists a closed Lipschitz curve in $(H, d_c)$ which does not bound an integral $2$-current in $(H, d_c)$, where $d_c$ is the Carnot-Carath\'eodory distance associated with the left-invariant Riemannian metric. The super-quadratic growth of $\FA(r)$ then follows from \thmref{main-thm-intcurr-isop} or \propref{prop-Carnot-isop}.

From \thmref{thm:inaccessible} and the fact that central products of finitely generated nilpotent groups of step $2$ have Dehn function bounded above by $n^2 \log n$, see \cite{Olshanskii-Sapir, Young-scaled-relators}, we obtain:

\bc\label{cor-lower-upper-inaccessible}
 Let $G_U$ and $H$ be as in \thmref{thm:inaccessible}, for some $m\geq 2$, and suppose the Lie algebra of $G_U$ has a basis with rational structure constants. Then the filling area and Dehn functions of $H$ satisfy
 \begin{equation}\label{eq:both-bounds}
  r^2 \varrho(r)\leq \FA(r) \leq \FA_0(r) \leq C r^2\log r
 \end{equation}
 for all $r\geq 2$, where $\varrho$ is a function satisfying $\varrho(r)\to\infty$ as $r\to\infty$.
\ec

Note that if $\mathfrak{g}$ and $U$ are as in \exref{example:inaccessible} for some $m\geq 2$ then, in particular, $G_U$ and its $m$-th central power $H= G_U\times_Z\dots\times_Z G_U$ satisfy the hypotheses of \corref{cor-lower-upper-inaccessible}. Consequently, the filling area and Dehn functions of $H$ satisfy \eqref{eq:both-bounds}.
This answers in the negative the question raised in \cite{Young-scaled-relators} whether the Carnot group $\Gamma = G_U\times_Z G_U$, where $\mathfrak{g}$ and $U$ are as in \exref{example:inaccessible} with $k=10$ and $m=2$, has quadratic Dehn function.

\begin{proof}[{Proof of \corref{cor-lower-upper-inaccessible}}]
 The lower bound for $\FA(r)$ in \eqref{eq:both-bounds} comes from \thmref{thm:inaccessible}. 
 In order to prove the upper bound for $\FA_0(r)$ let $\mathfrak{g}_U$ denote the Lie algebra of $G_U$ with grading $\mathfrak{g}_U = V_1\oplus V'_2$. Since $\mathfrak{g}_U$ has a basis with rational structure constants, there exists a basis of vectors in $V_1$ which generate a lattice $\Gamma$ in $G_U$. Now, there exists an injective homomorphism from the central product $\Gamma':= \Gamma \times_Z\dots\times_Z \Gamma$ of $m$ copies of $\Gamma$ to $H$ whose image is a lattice in $H$. Therefore, the Dehn function $\delta_{\Gamma'}(n)$ of $\Gamma'$ satisfies $$\delta_{\Gamma'}(n)\sim \FA_0(n).$$ Since $\Gamma'$ is a central product of $m\geq 2$ copies of a nilpotent group of step $2$, we have that $\delta_{\Gamma'}(n)\preceq n^2\log n$. This was first proved by Ol'shanskii-Sapir using the techniques of \cite{Olshanskii-Sapir}, see the remark on page 927 in \cite{Olshanskii-Sapir}. More recently, Young gave a proof of this in \cite{Young-scaled-relators} using different methods. It follows that there exists a suitable constant $C$ such that $\FA_0(r)$ satisfies
 \begin{equation*}
  \FA_0(r)\leq Cr^2\log r
 \end{equation*}
 for all $r\geq 2$. This concludes the proof.
\end{proof}

\thmref{thm:fin-gen-nilpotent} is a direct consequence of (the proof) of \corref{cor-lower-upper-inaccessible}. Indeed, for the finitely generated $\Gamma':= \Gamma \times_Z\dots\times_Z \Gamma$ appearing in the above proof, we have
 \begin{equation*}
  n^2 \varrho(n) \preceq \delta_\Gamma(n) \preceq n^2\log n
 \end{equation*}
 for some function $\varrho$ with $\varrho(n)\to\infty$ as $n\to\infty$.

\subsection{The proof of \thmref{thm:inaccessible}}

In the proof of the theorem it will often be useful to identify a given Carnot group $G$ via the exponential map with $(\mathfrak{g}, *)$, where $\mathfrak{g}$ is the Lie algebra of $G$ (viewed as a vector space) and $*$ is the multiplication on $\mathfrak{g}$ given by the Baker-Campbell-Hausdorff formula. The Lie algebra of $(\mathfrak{g}, *)$ is $\mathfrak{g}$ and the exponential map is simply the identity map on $\mathfrak{g}$.
If $G$ is of step $2$ then $*$ is given by
\begin{equation*}
 v * v'= v + v' + \frac{1}{2}[v, v'].
\end{equation*}

\begin{proof}[{Proof of \thmref{thm:inaccessible}}]
Let $G$, $\mathfrak{g} = V_1\oplus V_2$, $U$, $\mathfrak{g}_U$, $G_U$, and $H$ satisfy the hypotheses of the theorem. Suppose $U'\subset U$ is as in \defref{def:inaccessible} and let $U''\subset U$ be a subspace complementary to $U'$, thus $U = U' + U''$ and $U'\cap U'' = \{0\}$. 
Let $u_0\in U''$ with $u_0\not=0$ and choose $v_1, \dots, v_{2k}\in V_1$ with $$u_0= [v_1, v_2] + [v_3, v_4] + \dots + [v_{2k-1}, v_{2k}],$$
where $[\cdot, \cdot]$ is the Lie bracket on $\mathfrak{g}$. 
 Define, for every $j=1, \dots, k$, a piecewise affine curve $c_1^j:[0,4]\to V_1$ by
\begin{equation*}
 c_1^j(t) = \left\{\begin{array}{c@{\quad}l}
   t v_{2j-1} & 0\leq t<1\\
   v_{2j-1} + (t -1)v_{2j} & 1\leq t< 2\\
   (3-t)v_{2j-1} + v_{2j} & 2\leq t<3\\
   (4-t)v_{2j} & 3\leq t\leq 4,
 \end{array}\right.
\end{equation*}
and let $c_1:[0,1]\to V_1$ be the concatenation of the curves $c_1^j$, $j=1, \dots, k$, parametrized on $[0,1]$.  Clearly, $c_1$ is closed and piecewise affine. Define $\hat{c}_2: [0,1]\to V_2$ by 
\begin{equation*}
 \hat{c}_2(t) := \frac{1}{2} \int_0^{\, t}[c_1(s), \dot{c}_1(s)]ds
\end{equation*}
and observe that $\hat{c}_2(1) = u_0$.
Let now $$W_1:= V_1\oplus\dots\oplus V_1$$ be the direct sum of $m$ copies of $V_1$. Define quotient spaces $W_2:= V_2/U$ and $W'_2:= V_2/{U'}$ and denote by $P: V_2\to W_2$ and $P': V_2\to W'_2$ the natural projections. The Lie bracket  $[\cdot , \cdot]$ on $\mathfrak{g}$ gives rise to the bilinear map $W_1\times W_1\to V_2$, denoted by the same symbol $[\cdot , \cdot]$, defined by
\begin{equation}\label{eq:Lie-bracket-sum}
 [w, w'] = [v_1, v'_1] + \dots + [v_m, v'_m],
\end{equation}
where $w,w'\in W_1$ are of the form $w=v_1+\dots+ v_m$ and $w'=v'_1+\dots+ v'_m$ with $v_i$ and $v'_i$ in the $i$-th copy of $V_1$.
Clearly, $H$ has Lie algebra $\mathfrak{h} = W_1\oplus W_2$ with Lie bracket $$[w_1+ w_2, w'_1+w'_2]_{\mathfrak{h}}:= P ([w_1 , w'_1]),$$
where $w_i, w'_i\in W_i$ for $i=1,2$.
It follows from the $m$-inaccessibility of $U$ that for any $w_1, w'_1\in W_1$ we have 
\begin{equation}\label{eq:inacc-equiv}
[w_1, w'_1]_{\mathfrak{h}} = 0\quad\text{if and only if}\quad P'([w_1, w'_1]) = 0.
\end{equation}

Identifying $V_1$ with the first of the $m$ copies of $V_1$ in $W_1$, we may view $c_1$ as a curve in $W_1$. Define a curve $c:= [0,1]\to H$ by $c:= c_1 + P\circ\hat{c}_2$ and note that $c(1) = c_1(1) + P(u_0) = 0$.
Let $\|\cdot\|$ be a Euclidean norm on $W_1$ and let $d_0$ be the distance coming from a left-invariant Riemannian metric on $H$ which, restricted to $W_1\subset T_0H$, induces $\|\cdot\|$. Let furthermore $d_c$ be the Carnot-Carath\'eodory metric on $H$ associated with $\|\cdot\|$. Set $Y:= (H, d_c)$.
Clearly, $c$ is a closed Lipschitz curve in $Y$.  Set $T:= c_\#\Lbrack\chi_{[0,1]}\Rbrack$. We claim that there does not exist $S\in\intcurr_2(Y)$ with $\bdry S = T$. We argue by contradiction and assume such an $S\in\intcurr_2(Y)$ with $\bdry S = T$ exists. By \cite[Theorem 4.5]{Ambr-Kirch-curr}, $S$ is of the form
 \begin{equation*}
  S = \sum_{i=1}^\infty \varphi_{i\#}\Lbrack\theta_i\Rbrack
 \end{equation*}
 for some Lipschitz maps $\varphi_i: K_i\to Y$ with $K_i\subset\R^2$ compact and $\theta_i\in L^1(K_i, \Z)$. 
 Denote by $\pi: Y\to (W_1, \|\cdot\|)$ the projection given by $\pi(w):= w_1$ for 
 $w= w_1+w_2$ with $w_i\in W_i$, $i=1, 2$. As above, $H$ is identified with $W_1\oplus W_2$. It is straightforward to check that $\pi$ is $1$-Lipschitz and Pansu-differentiable at every point $y\in Y$ with Pansu-differential $d^P\pi(y): H\to W_1$ given by
 \begin{equation}\label{eq:Pansu-diff-pi}
  d^P\pi(x)(w) = w_1
 \end{equation}
 for $w = w_1 + w_2$ with $w_i\in W_i$, $i=1, 2$.
Set $\bar{S}:= \pi_\#S$ and note that $\bar{S}\in\intcurr_2((W_1, \|\cdot\|))$ with $$\bdry \bar{S} = \pi_\# T = c_{1\#}\Lbrack\chi_{[0,1]}\Rbrack;$$ furthermore
\begin{equation*}
 \bar{S} = \sum_{i=1}^\infty \psi_{i\#}\Lbrack\theta_i\Rbrack
\end{equation*}
where $\psi_i:= \pi\circ\varphi_i.$
Fix $i$. We claim that at almost every $x\in K_i$
\begin{equation}\label{eq:claim-bracket}
 [d_x\psi_i(u), d_x\psi_i(v)]_{\mathfrak{h}} = 0
\end{equation}
for all $u, v\in \R^2$, where $d_x\psi_i$ is the classical derivative (which exists almost everywhere by Rademacher's theorem). Indeed, by Pansu's Rademacher-type theorem \cite{Pansu-CC} and its generalization to Lipschitz maps defined only on measurable sets, see \cite{Vodop} and \cite{Magnani-area-formula}, the Pansu-differential $d^P\varphi_i(x): \R^2\to H$ of $\varphi_i$ at $x$ exists for almost every $x\in K_i$ and is a Lie group homomorphism, equivariant with respect to the dilations. Viewed as a map between Lie algebras, $d^P\varphi_i(x): \R^2\to\mathfrak{h}$ is a Lie algebra homomorphism and therefore
\begin{equation}\label{eq:Pansu-diff-Lie}
 \left[d^P\varphi_i(x)(u), d^P\varphi_i(x)(v)\right]_{\mathfrak{h}} = d^P\varphi_i(x)([u,v]_{\R^2}) = 0
\end{equation}
for all $u,v\in\R^2$. Since $\pi$ is Pansu-differentiable at $\varphi_i(x)$, the chain rule yields
\begin{equation*}
 d_x\psi_i = d^P\pi(\varphi_i(x)) \circ d^P\varphi_i(x).
\end{equation*}
This together with \eqref{eq:Pansu-diff-pi} and \eqref{eq:Pansu-diff-Lie} yields \eqref{eq:claim-bracket}, as claimed. It now follows from  \eqref{eq:claim-bracket} and \eqref{eq:inacc-equiv} that for almost every $x\in K_i$ and all $u, v\in \R^2$
\begin{equation}\label{eq:vanishing-bracket}
  P'([d_x\psi_i(u), d_x\psi_i(v)])= 0.
\end{equation}

Finally, let $\{\xi_1, \dots, \xi_n\}$ be a basis for $W_1$ and let $\pi^j: W_1\to \R$ be the corresponding coordinate functions, that is, $\pi^j(r_1\xi_1+\dots+r_n\xi_n) = r_j$, for $j=1, \dots, n$. Let furthermore $Q: W'_2\to\R$ be a linear functional and set $Q':= Q\circ P'$. Define functions $f_j: W_1\to \R$ by
\begin{equation*}
 f_j(x): = Q'([x, \xi_j]).
\end{equation*}
Clearly, the functions $f_j$ and $\pi^j$ are Lipschitz when $W_1$ is equipped with the  norm $\|\cdot\|$. Furthermore, $f_j$ is bounded on $\spt \bdry \bar{S}$. We calculate
\begin{equation*}
 \begin{split}
 \int_0^1 Q'([c_1(t), \dot{c}_1(t)]) dt & = \sum_{j=1}^n \bdry \bar{S}(f_j, \pi^j)\\
  & = \sum_{j=1}^n \bar{S}(1, f_j, \pi^j)\\
   &= \sum_{i=1}^\infty \sum_{j=1}^n \int_{K_i} \theta_i(x) \det\left(\nabla((f_j, \pi^j)\circ\psi_i)(x)\right) dx.
 \end{split}
\end{equation*}
Since, by an easy computation and \eqref{eq:vanishing-bracket}, 
\begin{equation*}
 \sum_{j=1}^n \det\left(\nabla((f_j, \pi^j)\circ\psi_i)(x)\right) = 2Q'( [d_x\psi_i(e_1), d_x\psi_i(e_2)] )= 0
\end{equation*}
for almost every $x\in K_i$, where $e_1$ and $e_2$ are the standard basis vectors of $\R^2$, we obtain
\begin{equation*}
  \int_0^1 Q'([c_1(t), \dot{c}_1(t)]) dt= 0.
\end{equation*}
Since $Q$ was arbitrary this shows that $P'(\hat{c}_2(1))= 0,$ a contradiction since $P'(\hat{c}_2(1)) = P'(u_0)\not=0$. This shows that there does not exist $S\in\intcurr_2(Y)$ with $\bdry S = T$. In particular, $Y$ does not admit a quadratic isoperimetric inequality for $\intcurr_1(Y)$. Since, by Pansu's result \thmref{pansu-thm}, the unique asymptotic cone of $X:=(H, d_0)$ is $Y$  it follows from  \thmref{main-thm-intcurr-isop} that $X$ does not admit a quadratic isoperimetric inequality for $\intcurr_1(X)$. This follows alternatively from \propref{prop-Carnot-isop}. Consequently, by Lemma 2.3 in \cite{Wenger-Gromov-hyp-isop}, the filling area function $\FA(r)$ of $X$ cannot be bounded by $Cr^2$ for any $C$. This completes the proof.
\end{proof}

\section{Another lower bound}\label{section-thm-subalgebra}

In this section we use similar arguments as above to prove super-quadratic lower bounds for the growth of the filling area functions for another class of Carnot groups. We prove:

\bt\label{thm-subalgebra}
 Let $G$ be a Carnot group of step $k$, endowed with a left-invariant Riemannian metric.  Let $$\mathfrak{g} = V_1 \oplus \dots \oplus V_k$$ be a grading of the Lie algebra $\mathfrak{g}$ of $G$. If $V_1$ does not contain a $2$-dimensional subalgebra then the filling area function $\FA(r)$ of $G$ grows strictly faster than quadratically:
 \begin{equation*}
  \frac{\FA(r)}{r^2}\to\infty\quad\text{as $r\to\infty$.}
 \end{equation*}
 \et

\begin{proof}
 By Pansu's result, \thmref{pansu-thm}, the unique asymptotic cone of $X:= (G, d_0)$ is $Y:=(G, d_c)$, where $d_c$ is the associated Carnot-Carath\'eodory distance. Since $Y$ is geodesic and is not a metric tree, it follows for example from Proposition 3.1 in \cite{Wenger-tree} that there exists $T\in\intcurr_1(Y)$ with $\bdry T=0$ and $T\not=0$. Suppose there exists $S\in\intcurr_2(Y)$ with $\bdry S = T$. Then $\|S\|$ is concentrated on a countably $\hm^2$-rectifiable subset $A\subset Y$ and is absolutely continuous with respect to $\hm^2$.  Since $Y$ is purely $2$-unrectifiable by \cite{Magnani} (see also \cite{AmKi-rect} for the case of the first Heisenberg group) it follows that $\hm^2(A)=0$ and thus also $S=0$. As a consequence, we obtain that $T = \bdry S = 0$, a contradiction. It thus follows that there exists no $S$ with $\bdry S = T$ and, in particular, $Y$ cannot admit a quadratic isoperimetric inequality for $\intcurr_1(Y)$. \thmref{main-thm-intcurr-isop} therefore shows that $X$ does not admit a quadratic isoperimetric inequality for $\intcurr_1(X)$. By Lemma 2.3 in \cite{Wenger-Gromov-hyp-isop}, the filling area function $\FA(r)$ of $X$ cannot be bounded by $Cr^2$ for any $C$.
\end{proof}

Simple examples of Carnot groups satisfying the hypotheses of \thmref{thm-subalgebra} are the first Heisenberg group and its generalizations using quaternions and octonions. It is known that the first Heisenberg group and the quaternionic Heisenberg group have cubic Dehn function, see \cite{Pittet-homogeneous-nilpotent} for the quaternionic case. In \cite{Pittet-homogeneous-nilpotent} it is also claimed that the octonionic Heisenberg group has cubic Dehn function. As was pointed out in \cite{Leuzinger-Pittet} there is a sign error in the proof of the octonionic case, so that the best previously known lower bound is quadratic. From \thmref{thm-subalgebra} we obtain a super-quadratic lower bound:

\bc
  The filling area function $\FA(r)$ of the octonionic Heisenberg group endowed with a left-invariant Riemannian metric grows strictly faster than quadratically.
\ec

For completeness, we recall the definition of the quaternionic and octonionic Heisenberg group. Denote by $\C$ the complex numbers, by $\Qu$ the quaternions, and by $\Oc$ the octonions. Recall that using the Cayley-Dickson construction, one obtains $\K_1:= \C$ from $\K_0:= \R$, $\K_2:= \Qu$ from $\K_1$, and $\K_3:= \Oc$ from $\K_2$ as follows: for $a\in\K_0$ define the conjugate $a^*$ of $a$ by $a^*:= a$ and the imaginary part of $a$ by $\Impart(a):=0$. Suppose that  for some $i\geq 1$, $\K_{i-1}$ has been defined together with multiplication, the conjugate $a^*$ of $a$, and the imaginary part $\Impart(a)$. Set $\K_i:= \{(a, b): a, b\in \K_{i-1}\}$ and define multiplication on $\K_i$ by
\begin{equation*}
 (a,b)(c,d):= (ac - db^*, a^*d+cb).
\end{equation*}
For $(a,b)\in\K_i$, define the conjugate by
\begin{equation*}
 (a,b)^*:= (a^*, -b)
\end{equation*}
and the imaginary part by $$\Impart(a,b):= (\Impart(a), b).$$

Now, set $\Li_i:= \Impart(\K_i)$ and define a stratified nilpotent Lie algebra  $\mathfrak{g}_i$ over $\R$ of step $2$ by
\begin{equation*}
 \mathfrak{g}_i:= \K_i\oplus \Li_i,
\end{equation*}
where the Lie bracket on $\mathfrak{g}_i$ is defined by
\begin{equation*}
 [z\oplus z', w\oplus w']:= \Impart(zw^*) \in \Li_i.
\end{equation*}

It is not difficult to check that for $z,w\in \K_i$ we have $\Impart(zw^*) = 0$ if and only if there exists $\lambda\in \R$ such that $w = \lambda z$. In particular, the first layer $\K_i$ of $\mathfrak{g}_i$ does not contain a $2$-dimensional subalgebra.

%\begin{proof}
% For $i=0$ it follows from the definition. Suppose now the lemma is true for some $i-1$ with $i\geq 1$. We prove it for $i$.
% The direction ``$\Longleftarrow$'' is clear. In order to prove ``$\Longrightarrow$'' let $z=(a,b)$ and $w=(c,d)$ with $a,b,c,d\in\K_{i-1}$ and compute 
% \begin{equation*}
%  zw^* = (ac^* + db^*, -a^*d+c^*b).
% \end{equation*}
%It follows that $\Impart(zw^*)=0$ if and only if
%\begin{equation*}
% a^*d = c^*b\quad\text{ and }\quad \Impart(ac^* + db^*) = 0.
%\end{equation*}
%Suppose first that $a=0$ and $b\not=0$. Then $\Impart(db^*) = 0$ and thus, since the lemma is assumed to be true for $i-1$, there exists $\lambda\in\R$ such that $d=\lambda b$. Since also $c^*b=0$ we have $c=0$. This proves that indeed $w=\lambda z$. Suppose now that $a\not=0$. It follows that
%\begin{equation*}
% d = \frac{1}{aa^*}ac^*b
%\end{equation*}
%and hence
%\begin{equation*}
% 0 = \Impart(ac^* + db^*) = \Impart\left(\left(1+\frac{bb^*}{aa^*}\right)ac^*\right) = \left(1+\frac{bb^*}{aa^*}\right)\Impart(ac^*).
%\end{equation*}
%From this it follows that $c=\lambda a$ for some $\lambda\in\R$ and hence $$d=\frac{1}{aa^*}ac^*b=\frac{aa^*}{aa^*}\lambda b = \lambda b.$$
%It follows again that $w=\lambda z$. This completes the proof.
%\end{proof}

%
\end{document}